\documentstyle[12pt]{article}
\title{On the Thomae formula for $Z_N$ curves}
\author{Atsushi Nakayashiki\\
Graduate School of Mathematics\\
Kyushu University}
\date{}

\begin{document}

\def\bcg{{\bf C}^g}
\def\bzg{{\bf Z}^g}
\def\ep{\epsilon}
\newcommand{\ch}[2]{
\left\{
\begin{array}{c}
#1\\ #2
\end{array}
\right\}_\tau}
\newcommand{\cht}[2]{
\left\{
\begin{array}{c}
#1\\ #2
\end{array}
\right\}_{\tau(t)}}
\newcommand{\th}[2]{
\theta\left[
\begin{array}{c}
#1 \\ #2
\end{array}
\right]
}
\newcommand{\bic}[2]{
\left(
\begin{array}{c}
#1 \\ #2
\end{array}
\right)
}
\def\la{\lambda}
\def\La{\Lambda}
\def\kp{\kappa}
\def\ra{\longrightarrow}
\def\ot{\otimes}
\newtheorem{prop}{Proposition}
\newtheorem{cor}{Corollary}
\newtheorem{lem}{Lemma}
\newtheorem{theo}{Theorem}
\newtheorem{definition}{Definition}
\newtheorem{conj}{Conjecture}

\maketitle

\begin{abstract}
We shall give an elementary and rigorous proof
of the Thomae formula for ${\bf Z}_N$ curves which
was discovered by Bershadsky and Radul \cite{BR1,BR2}.
Instead of using the determinant of the Laplacian
we use the traditional variational method
which goes back to Riemann, Thomae, Fuchs.
In the proof we made explicit the algebraic expression
of the chiral Szeg\"{o} kernels and proves the vanishing
of zero values of derivatives of 
theta functions with ${\bf Z}_N$ invariant $1/2N$
characteristics.
\end{abstract}
\par

\vskip4mm
\noindent
{\Large\bf 0 \hskip4mm Introduction}
\vskip4mm
\noindent
In \cite{BR1,BR2} Bershadsky and Radul discovered a
generalization of Thomae formula for ${\bf Z}_N$ curve
$s^N=f(z)=\prod_{i=1}^{Nm}(z-\la_i)$
(Theorem \ref{Thomae} in section 7 below).
The original Thomae formula is the case of hyperelliptic curves $N=2$
and takes the form
\begin{eqnarray}
&&
\theta[e](0)^8=\Big({\det A \over (2\pi i)^{m-1}}\Big)^4
\prod_{k<l}(\la_{i_k}-\la_{i_l})^2(\la_{j_k}-\la_{j_l})^2,
\nonumber
\end{eqnarray}
where $e$ is a non-singular even half period corresponding to the
partition of the branch points 
$\{1,\cdots,2m\}=\{i_1<\cdots<i_m\}\sqcup\{j_1<\cdots<j_m\}$,
$\{A_i,B_j\}$ a canonical homology basis and 
$A=(\int_{A_i}z^{j-1}/s)_{1\leq i,j\leq m-1}$.
This formula expresses the zero values of the Riemann theta 
functions with half characteristics as functions of
branch points.
Thomae formula was used to give generators of the affine ring of 
the moduli space
of hyperelliptic curves with level two structure \cite{M} 
in terms of theta constants, 
to give a generalization of the $\la$ function
of an elliptic curve \cite{M}(Umemura's appendix).
Beside those, F. Smirnov \cite{S} derived a beautiful theta formula
for the solution of the $sl_2$ Knizhnik-Zamolodchikov equation on level
zero using the Thomae formula.
For the generalized Thomae formula of ${\bf Z}_N$ curves
similar results are expected.
As for the generalization of $\lambda$ function for ${\bf Z}_N$
curves there are several results \cite{G,Far}
based on a different approach.
To study the generalization of Smirnov's formula to the case
of $sl_N$ is a main motivation for the present work.
It will be studied in a forthcomming paper.

Let us comment on the proof of the generalized Thomae formula.
Bershadsky and Radul evaluated, in two ways, 
the determinant of the Laplacian acting on 
some line bundle on a ${\bf Z}_N$ curve 
and compared them to obtain the formula.
However they used a path integral description of 
the correlation function of conformal fields to identify
it with the determinant of the Laplacian. 
Hence their proof does not seem mathematically rigorous.
It may be possible to make their proof rigorous
using the theory of determinants and Green functions only,
that is without the path integral.

Instead of going in that manner, here we shall give a rigorous and
elementary proof of the formula.
Our proof is based on the traditional variational method
which goes back to Riemann \cite{R}, 
Thomae \cite{T1,T2}, L. Fuchs \cite{F1,F2}.
The role of determinats and path integral is then replaced by
Fay's formula \cite{F} relating the Szeg\"{o} kernel and 
the canonical symmetric differential.
The strategy of the proof itself is similar to that of \cite{BR1,BR2}.

The particularity of our proof is to compare the analytic
and the algebraic expressions not only in the final formula
but also in each step of the proof.
As a corollary of those comparison the vanishing of the zero value
of the first order derivatives of theta functions with 
non-singular $1/2N$ characteristics is obatined.
This result is in turn used to prove the generalized Thomae formula.
Hence our proof clarifies some special aspects of theta functions
behind the generalized Thomae formula.
We also reveal a property of the proportionality constants
appeared in the Thomae formula for ${\bf Z}_N$ curves which
was not treated in \cite{BR1,BR2}.

\vskip2mm
Now the present paper is organized in the following manner.
In section 1 we gather necessary notation and formulas
concerning Riemann surfaces and theta functions following
the Fay's book \cite{F}.
The ${\bf Z}_N$ invariant $1/2N$ periods are studied
in section 2.
The algebraic expression for the chiral Szeg\"{o} kernel is
given in section 3.
Section 4 is devoted to the explanation of the canonical
differential and Fay's formula relating it with the 
chiral Szeg\"{o} kernel.
In section 5 the algebraic expression of the canonical differential
is studied.
The variation of the period matrix is studied in section 6.
In section 7 the generalized Thomae formula up to moduli
independent constant is proved.
The property of the proportionality constants is studied
in section 8.
In section 9 the examples of Thomae formula for small $N$'s
are given.
\vskip1cm

\section{Theta function}
\par
In this paper we mainly follow the notations of the
Fay's book \cite{F} which we summarize here.
Let $\tau$ be the $g$ by $g$ symmetric matrix whose real part is
negative definite. Any element $e\in \bcg$ is uniquely
expressed as 
\begin{eqnarray}
&&
e=\ch{\delta}{\ep}
=2\pi i\ep+\delta\tau,
\label{char}
\end{eqnarray}
with $\ep,\delta\in {\bf R}^g$. Here the vectors $\ep,\delta$ etc.
are all row vectors.
We call $\ep,\delta$ the characteristics of $e$.
The theta function with characteristics is defined by
\begin{eqnarray}
\th{\delta}{\ep}(z)&=&
\sum_{m\in\bzg}\exp(
{1\over2}(m+\delta)\tau(m+\delta)^{t}
+(z+2\pi i\ep)(m+\delta)^t
)
\nonumber
\\
&=&\exp(
{1\over2}\delta\tau\delta^t+
(z+2\pi i\ep)\delta^t
)\theta(z+e),
\nonumber
\end{eqnarray}
where
\begin{eqnarray}
&&
\th{0}{0}(z)=\theta(z),
\nonumber
\end{eqnarray}
and $e$ is determined by (\ref{char}).
We sometimes use $\theta[e](z)$ instead of 
$\th{\delta}{\ep}(z)$.
The transformation property is
\begin{eqnarray}
&&
\th{\delta}{\ep}(z+2\pi i\la+\kp\tau)
\nonumber
\\
&&
=
\exp(-{1\over2}\kp\tau\kp^t-z\kp^t+2\pi i(\delta\la^t-\ep\kp^t))
\th{\delta}{\ep}(z),
\nonumber
\end{eqnarray}
for $\la,\kp\in\bzg$.
We shall list some of the properties which easily follow from
the definitions:
\begin{eqnarray}
&&
\th{\delta+m}{\ep+n}(z)=
\exp(2\pi in\delta^t)\th{\delta}{\ep}(z),
\label{charprop1}
\\
&&
\th{-\delta}{-\ep}(0)=
\th{\delta}{\ep}(0),
\nonumber
\\
&&
{\th{\delta}{\ep}(z)\over \th{\delta}{\ep}(0)}
=
{\th{\delta+m}{\ep+n}(z)\over \th{\delta+m}{\ep+n}(0)},
\nonumber
\end{eqnarray}
for $m,n\in\bzg$.

\par

Let $C$ be a compact Riemann surface of genus $g$.
Let us fix a marking of $C$ \cite{Gun}(\S1). 
That means, we fix a canonical basis $\{A_i,B_j\}$ of $\pi_1(C)$,
a base point $P_0\in C$ and a base point in the universal
cover $\tilde{C}$ which lies over $P_0$.
We assume that the tails of $A_i,B_j$ are joined to $P_0$.
Then we can canonically identify the covering transformation group and
the fundamental group $\pi_1(C,P_0)$.
We also identify holomorphic $1$-forms on $C$ with holomorphic $1$-forms
on $\tilde{C}$ invariant under the action of $\pi_1(C)$.
Let us denote by $\pi:\tilde{C}\ra C$ the projection and 
by $J(C)$ the Jacobian variety of $C$ which is the set of linear
equivalence classes of degree $0$ divisors on $C$.
In the following sections we always assume one marking of $C$.

Let $\{v_j\}$  be the basis of the normalized holomorphic 1-forms.
The normalization is
\begin{eqnarray}
&&
\int_{A_j}v_k=2\pi i\delta_{jk},
\nonumber
\end{eqnarray}
and set
\begin{eqnarray}
&&\int_{B_j}v_k=\tau_{jk}.
\nonumber
\end{eqnarray}
A flat line bundle on $C$ is described by the character
of the fundamental group $\chi:\pi_1(C)\ra{\bf C}^\ast$,
where ${\bf C}^\ast$ is the multiplicative group of
non-zero complex numbers.
The two representation $\chi_1$ and $\chi_2$ defines
a holomorphically equivalent line bundle if and only if
there exists an holomorphic 1-form $\omega$ such that
\begin{eqnarray}
&&
\chi_1(\gamma)\chi_2(\gamma)^{-1}
=\exp(\int_\gamma\omega)
\nonumber
\end{eqnarray}
for any $\gamma\in\pi_1(C)$.
Let ${\cal A}$ and ${\cal B}$ be positive divisors of the
same degree, say $d$, and set ${\cal A}=\sum_{i=1}^d P_i$,
${\cal B}=\sum_{i=1}^d R_i$.
Let us fix points $\tilde{P}_i,\tilde{R}_j$ in $\tilde{C}$
so that they lie over $P_i,R_j$. 
Let us set
\begin{eqnarray}
&&
\int_{{\cal A}}^{{\cal B}}v_i=
\sum_{j=1}^d\int_{\tilde{P}_j}^{\tilde{R}_j}v_i,
\nonumber
\end{eqnarray}
where the integration in the right hand side is taken in $\tilde{C}$.
Then the flat line bundle corresponding to the degree $0$
divisor ${\cal B}-{\cal A}$ is described by
\begin{eqnarray}
&&
\chi(A_i)=1,
\qquad
\chi(B_i)=\exp(\int_{{\cal A}}^{{\cal B}}v_i).
\nonumber
\end{eqnarray}
Another choice of $\tilde{P}_i,\tilde{R}_j$ gives 
an equivalent line bundle.

We say $\chi$ is unitary if the image is contained in
the unitary group $U(1)$.
The following proposition is well known and easily proved.

\begin{prop}
For an isomorphism class of flat line bundles there
exists a unique unitary representation $\chi$ which defines
the line bundle belonging to that class.
\end{prop}
\vskip2mm

If we take $\delta,\epsilon\in{\bf R}^g$ such that
\begin{eqnarray}
&&
(\int_{{\cal A}}^{{\cal B}}v_1,\cdots,\int_{{\cal A}}^{{\cal B}}v_g)
=\ch{\delta}{\ep}
\nonumber
\end{eqnarray}
as a point on the Jacobian variety of $C$,
then the corresponding unitary representation $\tilde{\chi}$ 
is given by
\begin{eqnarray}
&&
\tilde{\chi}(A_j)=\exp(-2\pi i\delta_j),
\qquad
\tilde{\chi}(B_j)=\exp(2\pi i\ep_j).
\label{unitary}
\end{eqnarray}
The multiplicative meromorphic function described by $\tilde{\chi}$
is, for example, given by
\begin{eqnarray}
&&
{
\th{-\delta}{-\ep}(\int_{z_0}^xv-\alpha)\over
\theta(\int_{z_0}^xv-\alpha)
},
\nonumber
\end{eqnarray}
where $v$ is the vector of normalized holomorphic 1-forms,
$x\in \tilde{C}$, $\alpha\in \bcg$ and the integration path
is taken in $\tilde{C}$.

We denote by $\Delta$ the Riemann divisor for our choice
of the canonical homology basis which satisfies
\begin{eqnarray}
&&
2\Delta\equiv K_C.
\nonumber
\end{eqnarray}
Here $K_C$ is the divisor class of the canonical bundle
of $C$ and $\equiv$ means the linear equivalence.
Let $L_0$ be the degree $g-1$ line bundle corresponding to $\Delta$.
For a divisor $\alpha$ with degree $0$ let us denote by
${\cal L}_\alpha$ the corresponding flat line bundle and set
$L_\alpha={\cal L}_\alpha\ot L_0$.
For a non-singular odd half period $\alpha$ let $h_\alpha$ be
the section of $L_\alpha$ which satisfies
\begin{eqnarray}
&&
h_\alpha^2(x)=
\sum_{j=1}^g
{\partial\theta[\alpha]\over \partial z_j}(0)v_j(x).
\nonumber
\end{eqnarray}
Then the prime form is defined by
\begin{eqnarray}
E(x,y) &=& {\theta[\alpha](y-x)\over h_\alpha(x)h_\alpha(y)},
\nonumber
\\
y-x &=& \int_x^y v,
\nonumber
\end{eqnarray}
where $x,y\in \tilde{C}$ and $v=(v_1(x),\cdots,v_g(x))$.
Let $\pi_j$ be the projection from $C\times C$ to the $j$-th
component and $\delta:C\times C\ra J(C)$ the map $(x,y)\mapsto y-x$.
Then $E(x,y)$ can be considered as a section of the line bundle
$\pi_1^\ast L_\alpha\ot\pi_2^\ast L_\alpha\ot\delta^\ast\Theta$,
where $\Theta$ is the line bundle on $J(C)$ defined by the
theta divisor.
Let us fix the transformation property of the half differential
on $\tilde{C}$ under the action of $\pi_1(C)$ so that the
section of $\pi^\ast L_0$ is invariant.
This means, in particular, that $E(x,y)$ transforms under
the action of $A_i$, $B_i$ in $y$ as
\begin{eqnarray}
&&
E(x,y+A_i)=E(x,y),
\quad
E(x,y+B_i)=\exp(-{\tau_{ii}\over2}-\int_x^y v_i)E(x,y).
\nonumber
\end{eqnarray}
Here we denote the action of $A_i$, $B_i$ in an additive manner.
The prime form has the nice expansion as follows.
Let $u$ be a local coordinate around $P\in \tilde{C}$.
Then the expansion of $E(x,y)$ in $u(y)$ at $u(x)$ takes the form
\begin{eqnarray}
&&
E(x,y)\sqrt{du(x)}\sqrt{du(x)}
=u(y)-u(x)+O\big((u(y)-u(x))^3\big).
\label{primeexp}
\end{eqnarray}
Since the expansion is of local nature we
sometimes use the way of saying 
that $P\in C$, the local corrdinate $u$ around $P$
and the expansion in $u(y)$ at $u(x)$ etc.
\vskip1cm

\section{${\bf Z}_N$ curve and ${1\over2N}$ period}
\par
Let us consider the plane algebraic curve 
$s^N=f(z)=\prod_{j=1}^{Nm}(z-\la_i)$.
We compactify it by adding $N$ infinity point 
$\infty^{(1)},\cdots,\infty^{(N)}$ and denote 
the compact Riemann surface by $C$.
The genus $g$ of $C$ is $g=1/2(N-1)(Nm-2)$.
The $N$-cyclic automorphism $\phi$ of $C$ is defined by
$\phi:(z,s)\mapsto(z,\omega s)$, where $\omega$ is the $N$-th
primitive root of unity.
There are $Nm$ branch points $Q_1,\cdots,Q_{Nm}$ whose projection
to $z$ coordinate are $\la_1,\cdots,\la_{Nm}$.

The basis of holomorphic 1-forms on $C$ is given by
\begin{eqnarray}
&&
w^{(\alpha)}_\beta={z^{\beta-1}dz\over s^{\alpha}}
\quad
1\leq\alpha\leq N-1,
\quad
1\leq\beta\leq \alpha m-1.
\nonumber
\end{eqnarray}

Let us describe the divisors which we need and their relations.
The following lemma is easily proved.
\begin{lem}
For any $P\in C$ the linear equivalence 
class $P+\phi(P)+\cdots+\phi^{N-1}(P)$
does not depend on the point $P$.
\end{lem}

We set 
\begin{eqnarray}
&&
D\equiv P+\phi(P)+\cdots+\phi^{N-1}(P).
\nonumber
\end{eqnarray}

The following lemma is easily proved.

\begin{lem}\label{divrels}
The following relations hold.
\begin{description}
\item[1.] $D\equiv NQ_i\equiv \infty^{(1)}+\cdots+\infty^{(N)}$ 
for any $i$.
\item[2.] $K_C\equiv (L-1)D$, where $L=(N-1)m-1$.
\item[3.] $\sum_{j=1}^{Nm}Q_j\equiv mD$.
\end{description}
\end{lem}
\vskip2mm

Following Bershadsky-Radul\cite{BR2} we shall describe
the important object of our study, the ${\bf Z}_N$ invariant
$1/N$ or $1/2N$ periods.
Let us consider an ordered partition $\La=(\La_0,\cdots,\La_{N-1})$ of
$\{1,2,\cdots, Nm\}$
such that the number $|\La_i|$ of elements of $\La_i$
is equal to $m$ for any $i$.
With each $\La$ we associate the divisor class $e_\La$ by
\begin{eqnarray}
&&
e_\La \equiv \La_1+2\La_2+\cdots+(N-1)\La_{N-1}-D-\Delta,
\nonumber
\end{eqnarray}
where for a subset $S$ of $\{1,2,\cdots, Nm\}$
we set
\begin{eqnarray}
&&
S=\sum_{j\in S}Q_j.
\nonumber
\end{eqnarray} 
For a given $\La$ we denote by $\La(j)$ the ordered partition
\begin{eqnarray}
&&
\La(j)=(\La_j,\La_{j+1},\cdots,\La_{j-1}).
\nonumber
\end{eqnarray}
Here we consider the index of $\La_j$ by modulo $N$.
Then

\begin{prop}
For any ordered partition $\La$ we have
\begin{description}
\item[1.] $Ne_\La\equiv0$ for $N$ being even and $2Ne_\La\equiv0$
for $N$ being odd.
\item[2.] $e_\La\equiv e_{\La(2)}\equiv\cdots\equiv e_{\La(N)}$.
\item[3.] $-e_\La\equiv \La_{N-1}+2\La_{N-2}+\cdots+(N-1)\La_1-D-\Delta$.
\end{description}
\end{prop}
\vskip2mm

This proposition is easily proved using Lemma \ref{divrels}.
For $\La=(\La_0,\ldots,\La_{N-1})$ we set
\begin{eqnarray}
\La^{-}=(\La^{-}_0,\ldots,\La^{-}_{N-1})=(\La_0,\La_{N-1},\ldots,\La_1),
\nonumber
\end{eqnarray}
and $\La=\La^{+}$.
Let $\theta(z)$ be the theta function associated with
our choice of canonical homology basis.
Then

\begin{prop}\label{nonsing}
The $1/2N$ period $e_\La$ is non-singular, that means
\begin{eqnarray}
&&
\theta(e_\La)\neq0.
\nonumber
\end{eqnarray}
\end{prop}

This proposition was proved in \cite{BR2}.
One can find another proof in \cite{G} which is similar to
that of \cite{M} in the hyperelliptic case.
\vskip1cm

\section{Chiral Szeg\"{o} kernel}
\par

\newtheorem{defi}{Definition}
\begin{defi}
For $e\in{\bf C}^g$ satisfying $\theta(e)\neq0$ the chiral Szeg\"{o} kernel
$R(x,y\vert e)$ is defined by
\begin{eqnarray}
&&
R(x,y\vert e)=
{
\theta[e](y-x)\over\theta[e](0)E(x,y)
}
\qquad x,y\in \tilde{C}.
\nonumber
\end{eqnarray}
\end{defi}
\vskip2mm

We remark that $R(x,y\vert e)$ depends only on
the image of $e$ to the Jacobian variety $J(C)$.

We shall give an algebraic expression for $R(x,y\vert e_\La)$.
Let us set
\begin{eqnarray}
{\cal L}&=&\{-{N-1\over2},-{N-1\over2}+1,\cdots,{N-1\over2}\},
\label{calL}
\\
q_l(i)&=&{1-N\over 2N}+\Big\{{l+i+{N-1\over2}\over N}\Big\},
\label{qli}
\end{eqnarray}
for $l\in{\cal L}$ and $i\in{\bf Z}$.
Here $\{a\}=a-[a]$ is the fractional part of $a$, $[a]$ being the Gauss
symbol.
For an ordered partition $\La=(\La_0,\cdots,\La_{N-1})$ we
define the number $k_i, i=1,\cdots,Nm$ by
\begin{eqnarray}
&&
i\in\La_j\quad\hbox{if and only if }k_i=j.
\label{weight}
\end{eqnarray}
For each $l\in{\cal L}$ we set
\begin{eqnarray}
&&
f_l(x,\La)=\prod_{i=1}^{Nm}(z(x)-\la_i)^{q_l(k_i)}\sqrt{dz(x)}.
\nonumber
\end{eqnarray}
The following proposition was found in \cite{BR2}.

\begin{prop}\label{bundle}
$f_l(x,\La)$ is a meromorphic section of 
$L_{e_\La}$ whose divisor is
\begin{eqnarray}
&&
{\sf div} f_l=\La_{1-j}+2\La_{2-j}+\cdots+(N-1)\La_{-1-j}-
\sum_{k=1}^N\infty^{(k)},
\label{divisor}
\end{eqnarray}
where $l=-(N-1)/2+j$.
\end{prop}
\vskip2mm

Note that the chiral Szeg\"{o} kernel $R(x,y\vert e_\La)$ 
can be considered as a section of the line bundle 
$\pi_1^\ast L_{e_\La}\ot \pi_2^\ast L_{-e_\La}$,
where $\pi_i$ is the projection to the $i$-th component of 
$C\times C$.
Let us set
\begin{eqnarray}
&&
F(x,y\vert \La)={1\over N}
{\sum_{l\in{\cal L}}f_l(x,\La)f_{-l}(y,\La^{-})
\over z(y)-z(x)}.
\nonumber
\end{eqnarray}
Here the choice of the branch of $f_l(x,\La)$
should be specified as in (\ref{branch1}) and (\ref{branch2}).
Note that $F(x,y\vert \La)$ and 
$R(x,y\vert e_\La)$ can be considered as the sections of 
the same line bundle.
Then we have

\begin{theo}\label{mainth1}
For an ordered partition $\La$ we have
\begin{eqnarray}
&&
R(x,y\vert e_\La)=F(x,y\vert \La).
\nonumber
\end{eqnarray}
\end{theo}
\vskip2mm

As a corollary of this theorem we have the vanishing
of the theta derivative constants.

\begin{cor}\label{vanish}
For any ordered partition $\La$ we have
\begin{eqnarray}
&&
{\partial\theta[e_\La]\over\partial z_i}(0)=0
\qquad\hbox{for any $i$}.
\label{vanishing}
\end{eqnarray}
\end{cor}

Note that whether the right hand side of 
 (\ref{vanishing}) vanishes or not depends only on the divisor class
of $e_\La$ by (\ref{charprop1}).
Hence the statement has unambiguously a sense.
This curious result is a natural generalization of the
hyperelliptic case where $e_\La$ is a non-singular
even half period and the corollary is obvious.
For general $N$ we do not know whether $\theta[e_\La](z)$
is an even function.

\begin{lem}\label{exponent}
The following properties hold.
\begin{description}
\item[1.] $q_l(i)=q_{l^\prime}(i^\prime)$ if $i+l=i^\prime+l^\prime$.
\item[2.] $q_l(i+N)=q_l(i)$ for any $i$.
\item[3.] $\sum_{i=0}^{N-1}q_l(i)=0$.
\end{description}
\end{lem}
\vskip2mm
\noindent
Proof. The properties 1 and 2 are obvious.
Let us prove 3. Using 1 and 2 we have
\begin{eqnarray}
&&
\sum_{i=0}^{N-1}q_l(i)=\sum_{i=0}^{N-1}q_{-{N-1\over2}}(i)
=\sum_{i=0}^{N-1}\big({N-1\over2N}+{i\over2}\big)=0.
\nonumber
\end{eqnarray}
$\Box$
\vskip5mm

\noindent
Proof of Proposition \ref{bundle}. The meromorphy at points except
the branch points and $\infty^{(k)}$ is obvious.
We can take $t=(z-\la_i)^{1/N}$ as a local coordinate around $Q_i$.
Then
\begin{eqnarray}
&&
(z-\la_i)^{q_l(k_i)}\sqrt{dz}=t^{{N-1\over2}+Nq_l(k_i)}\sqrt{Ndt}.
\nonumber
\end{eqnarray}
If we write $l=-(N-1/2)+j$ $(0\leq j\leq N-1)$, we have
\begin{eqnarray}
&&
{N-1\over2}+Nq_l(k_i)=N\{{k_i+j\over N}\}.
\nonumber
\end{eqnarray}
At $\infty^{(k)}$ we can take $t=1/z$ as a local coordinate
and we have
\begin{eqnarray}
&&
f_l(x,\La)={1\over t}\sqrt{dt}(1+O(t))
\nonumber
\end{eqnarray}
by the property 3 of Lemma \ref{exponent}.
Hence $f_l$ is locally meromorphic on $C$
with the divisor (\ref{divisor}).
Let us consider $f_l(x,\La)^2$.
This is a multi-valued meromorphic 1-form with the divisor
\begin{eqnarray}
&&
2\big(\sum_{k=1}^{N-1}k\La_{k-j}
-\sum_{k=1}^N\infty^{(k)}\big)\equiv
2e_{\La(1-j)}+2\Delta\equiv 2e_\La+K_C.
\nonumber
\end{eqnarray}
Hence
\begin{eqnarray}
&&
f_l(x,\La)^2\in
H^0\big(C,
{\cal L}_{e_\La}^{\ot 2}\ot 
\Omega^1_C(-2\sum_{k=1}^{N-1}k\La_{k-j}+
2\sum_{k=1}^N\infty^{(k)})
\big).
\nonumber
\end{eqnarray}
Since
\begin{eqnarray}
&&
{\cal L}_{e_\La}^{\ot 2}\ot 
\Omega^1_C(-2\sum_{k=1}^{N-1}k\La_{k-j}+
2\sum_{k=1}^N\infty^{(k)})\simeq {\cal O}_C,
\nonumber
\end{eqnarray}
we have
\begin{eqnarray}
&&
H^0\big(C,
{\cal L}_{e_\La}^{\ot 2}\ot 
\Omega^1_C(-2\sum_{k=1}^{N-1}k\La_{k-j}+
2\sum_{k=1}^N\infty^{(k)})
\big)={\bf C}f_l(x,\La)^2.
\nonumber
\end{eqnarray}
Note that
\begin{eqnarray}
&&
H^0\big(C,
L_{e_\La}(-\sum_{k=1}^{N-1}k\La_{k-j}+
\sum_{k=1}^N\infty^{(k)})
\big)
\label{cohomology}
\end{eqnarray}
is one dimensional. 
Hence $f_l(x,\La)$ can be considered as an element of (\ref{cohomology}).
Thus Proposition \ref{bundle} is proved.
$\Box$

\begin{lem}
The following expression holds :
\begin{eqnarray}
&&
f_{-l}(y,\La^{-})=\prod_{i=1}^{Nm}(z(y)-\la_i)^{-q_l(k_i)}\sqrt{dz(y)}.
\nonumber
\end{eqnarray}
\end{lem}
\vskip2mm

\noindent
Proof. Recall that 
\begin{eqnarray}
&&
\La^{-}=(\La^{-}_0,\cdots,\La^{-}_{N-1}),
\quad
\La^{-}_j=\La_{N-j}.
\nonumber
\end{eqnarray}
Then we have
\begin{eqnarray}
&&
f_{-l}(y,\La^{-})=
\prod_{i=1}^{Nm}(z(y)-\la_i)^{q_{-l}(N-k_i)}\sqrt{dz(y)}.
\nonumber
\end{eqnarray}
Hence it is sufficient to prove
\begin{eqnarray}
&&
q_{-l}(N-i)=-q_l(i),
\nonumber
\end{eqnarray}
for any $l$ and $i$.
This can be easily checked. $\Box$

\begin{lem}
$F(x,y\vert \La)$ is regular outside the diagonal set $\{x=y\}$.
\end{lem}
\vskip2mm

\noindent
Proof. A priori we know that $F(x,y\vert \La)$ has poles
at most at $z(x)=z(y)$.
Hence it is sufficient to prove that $F(x,y\vert \La)$
is regular at $z(x)=z(y)$ and $x\neq y$.
If we write $l=-(N-1)/2+j$ we have
\begin{eqnarray}
q_l(k_i)-q_{-{N-1\over2}}(k_i)=
q_{-{N-1\over2}}(j+k_i)-q_{-{N-1\over2}}(k_i)
={j\over N}\hbox{ mod. ${\bf Z}$}.
\nonumber
\end{eqnarray}
Therefore we can set
\begin{eqnarray}
&&
q_l(k_i)-q_{-{N-1\over2}}(k_i)={j\over N}+r_{ij},
\quad
r_{ij}\in{\bf Z}.
\nonumber
\end{eqnarray}
Let us choose the branch of $f_l(x,\La)$, $f_l(x,\La^{-})$
such that the following equations hold:
\begin{eqnarray}
f_l(x,\La)&=&
(z(x)-\la_i)^{r_{ij}}s(x)^j
f_{-{N-1\over2}}(x,\La),
\label{branch1}
\\
f_{-l}(x,\La^{-})&=&
(z(x)-\la_i)^{-r_{ij}}s(x)^{-j}
f_{-{N-1\over2}}(x,\La^{-}).
\label{branch2}
\end{eqnarray}
Then we have
\begin{eqnarray}
&&
F(x,y\vert \La)
\nonumber
\\
&&=
f_{-{N-1\over2}}(x,\La)
f_{-{N-1\over2}}(y,\La^{-})
\sum_{j=0}^{N-1}\prod_{i=1}^{Nm}
\Big(
{z(x)-\la_i\over z(y)-\la_i}
\Big)^{r_{ij}}
\Big({s(x)\over s(y)}\Big)^j.
\nonumber
\end{eqnarray}
Now in the limit 
\begin{eqnarray}
&&
z(x)\ra z(y),
\quad
s(x)\ra \omega^r s(y),
\quad
1\leq r\leq N-1,
\nonumber
\end{eqnarray}
with $\omega=\exp{2\pi i/N}$ we have
\begin{eqnarray}
&&
F(x,y\vert \La)\ra
f_{-{N-1\over2}}(y^{(r)},\La)
f_{-{N-1\over2}}(y,\La^{-})
\sum_{j=0}^{N-1}\omega^r=0,
\nonumber
\end{eqnarray}
where $y^{(r)}=(z(y),\omega^rs(y))$. $\Box$

The following lemma is proved by a direct calculation.
\begin{lem}\label{expand1}
Let $P\in C$ be a non-branch point.
We can take $z$ to be a local coordinate around $p$.
Then the expansion in $z(y)$ at $z(x)$ takes the form
\begin{eqnarray}
&&
F(x,y\vert \La^\pm)
\nonumber
\\
&&=
{\sqrt{dz(x)}\sqrt{dz(y)}\over z(y)-z(x)}
\Big[1+
{1\over 2N}\sum_{i,j=1}^{Nm}
{q(k_i,k_j)\over(z(x)-\la_i)(z(x)-\la_j)}
(z(y)-z(x))^2+\cdots\Big],
\nonumber
\end{eqnarray}
where $q(i,j)=\sum_{l\in{\cal L}}q_l(i)q_l(j)$.
\end{lem}
\vskip2mm

The following lemma is a consequence of the expansion
(\ref{primeexp}) of the prime form $E(x,y)$.
\begin{lem}\label{expand3}
Under the same conditions of Lemma \ref{expand1} we have
\begin{eqnarray}
&&
R(x,y\vert e_\La)
\nonumber
\\
&&=
{\sqrt{dz(x)}\sqrt{dz(y)}\over z(y)-z(x)}
\Big[1+
\sum_{i=1}^{g}
{\partial\log\theta[e_\La]\over\partial z_i}(0)v_i(x)
(z(y)-z(x))+\cdots\Big],
\nonumber
\end{eqnarray}
where $v_i(x)$ means the coefficient of $dz(x)$ in $v_i(x)$.
\end{lem}
\vskip2mm

\noindent
Proof of Theorem \ref{mainth1}. 
Let $\chi$ be the unitary representation corresponding to 
${\cal L}_{e_\La}$. If we write
\begin{eqnarray}
&&
e_{\La}=\ch{\delta}{\ep},
\nonumber
\end{eqnarray}
then $\chi(A_i)$ and $\chi(B_j)$ are given by (\ref{unitary})
in section 1.
The transformation property of $R(x,y\vert e_\La)$ is
\begin{eqnarray}
&&
R(x+\gamma_1,y+\gamma_2 \vert e_\La)
=\chi(\gamma_1)\chi(\gamma_2)^{-1}R(x,y\vert e_\La),
\nonumber
\end{eqnarray}
for $\gamma_1,\gamma_2\in\pi_1(C)$.
On the other hand if we pull $F(x,y\vert \La)$ 
back to $\tilde{C}\times\tilde{C}$ then
\begin{eqnarray}
&&
F(x+\gamma_1,y+\gamma_2 \vert \La)
=\chi_1(\gamma_1)\chi_2(\gamma_2)F(x,y\vert \La),
\nonumber
\end{eqnarray}
for some unitary representation $\chi_1$ and $\chi_2$.
In fact if $x$ rounds a cycle of $C$
$f_l(x,\La)$ is multiplied by an appropriate
$2N$ the root of unity. The same is true for $y$.

Let us set
\begin{eqnarray}
\tilde{F}(x,y\vert \La)=
{R(x,y\vert e_\La)\over F(x,y\vert \La)}.
\nonumber
\end{eqnarray}
Then $\tilde{F}(x,y\vert \La)$
is the section of the trivial line bundle
and obeys the tensor product of unitary representations
of $\pi_1(C)\times\pi_1(C)$.
Hence $\tilde{F}(x,y\vert \La)$ is invariant 
under the action of $\pi_1(C)\times\pi_1(C)$.
This means that $R(x,y\vert e_\La)$ and $F(x,y\vert \La)$
have the same transformation property.
Therefore the function
\begin{eqnarray}
&&
I(x,y)=F(x,y\vert \La)-R(x,y\vert e_\La)
\nonumber
\end{eqnarray}
can be considered as a section of the line bundle
$
\pi_1^\ast L_{e_\La}\ot 
\pi_2^\ast L_{-e_\La}.
$
By Lemma \ref{expand1} and \ref{expand3} we know that $I(x,y)$
is holomorphic except $\cup_{i=1}^{Nm}\{Q_i\}\times \{Q_i\}$.
Since $I(x,y)$ is meromorphic on $C\times C$,
$I(x,y)$ has no singularity.
By Proposition \ref{nonsing},
$e_\La$ is non-singular which means
\begin{eqnarray}
&&
H^0(C,L_{e_\La})=0.
\nonumber
\end{eqnarray}
Hence
\begin{eqnarray}
&&
H^0(C\times C,\pi_1^\ast L_{e_\La}
\ot \pi_2^\ast L_{-e_\La})
=\pi_1^\ast H^0(C,L_{e_\La})\ot
\pi_2^\ast H^0(C,L_{-e_\La})=0.
\nonumber
\end{eqnarray}
Thus $I(x,y)=0$. $\Box$
\vskip5mm

\noindent
Proof of Corollary \ref{vanish}. This is a direct consequence
of Theorem \ref{mainth1}, Lemma \ref{expand1} and \ref{expand3}. $\Box$

\begin{cor}\label{expand2}
Under the same conditions and notations as in Lemma \ref{expand1},
then we have
\begin{eqnarray}
&&
R(x,y\vert e_\La)R(x,y\vert -e_\La)
\nonumber
\\
&&
={dz(x)dz(y)\over (z(y)-z(x))^2}
\Big[1+
{1\over N}\sum_{i,j=1}^{Nm}
{q(k_i,k_j)\over(z(x)-\la_i)(z(x)-\la_j)}
(z(y)-z(x))^2+\cdots\Big],
\nonumber
\end{eqnarray}
\end{cor}
\vskip1cm

\section{Canonical symmetric differential}
\par
The canonical symmetric differential $\omega(x,y)$ 
is defined by the following properties.
\vskip3mm
\begin{description}
\item[1] $\omega(x,y)$ is a meromorphic section
of $\pi_1^\ast\Omega_C^1\ot \pi_2^\ast\Omega_C^1$ on $C\times C$, where
$\pi_i$ is the projection to the $i$-th component of $C\times C$.
\item[2] $\omega(x,y)$ is holomorphic except the diagonal set $\{x=y\}$ where
it has a double pole. For $p\in C$ if we take a local coordinate $u$
around $p$ then the expansion in $u(x)$ at $u(y)$ takes the form
\begin{eqnarray}
&&
\omega(x,y)=\Big({1\over(u(x)-u(y))^2}+\hbox{regular}\Big)
du(x)du(y).
\nonumber
\end{eqnarray}
\item[3] The $A$ period in $x$ variable is zero:
\begin{eqnarray}
&&
\int_{A_j}\omega(x,y)=0
\quad \hbox{for any $j$}.
\nonumber
\end{eqnarray}
\item[4] $\omega(x,y)=\omega(y,x)$.
\end{description}
\vskip5mm

The following proposition is well known.

\begin{prop}
The canonical differential exists and is unique.
\end{prop}

In fact there is an analytical description of $\omega(x,y)$ in terms
of the theta function (see for example \cite{F} p26, Corollary 2.6) :
\begin{eqnarray}
&&
\omega(x,y)=d_xd_y\log E(x,y)
=-\sum_{i,j=1}^g{\partial^2\log\theta\over\partial z_i\partial z_j}
(y-x-f)v_i(x)v_j(y),
\nonumber
\end{eqnarray}
for any non-singular point $f\in(\Theta)$, 
where $(\Theta)=(\theta(z)=0)$.
The uniqueness can be easily proved using
$H^0(C\times C,\pi_1^\ast\Omega_C^1\ot \pi_2^\ast\Omega_C^1)
=\pi_1^\ast H^0(C,\Omega_C^1)\otimes \pi_2^\ast H^0(C,\Omega_C^1)$.

There is a remarkable identity due to Fay \cite{F} ( Corollary 2.12 )
connecting the chiral Szeg\"{o} kernel and the canonical 
symmetric differential.
The formula is
\begin{eqnarray}
&&
R(x,y\vert e)R(x,y\vert -e)
=
\omega(x,y)+
\sum_{i,j=1}^g
{\partial^2\log\theta[e]\over\partial z_i\partial z_j}
(0)v_i(x)v_j(y),
\label{remarkable}
\end{eqnarray}
for any $e\in {\bf C}^g$ such that $\theta(e)\neq0$.

For a non-branch point $P\in C$ we can take $z$ as
a local coordinate around $P$. Let us define
\begin{eqnarray}
&&G_z(z)=\lim_{y\rightarrow x}
\Big[
\omega(x,y)-{dz(x)dz(y)\over((z(y)-z(x))^2}
\Big].
\nonumber
\end{eqnarray}
It is known\cite{F,HS} that $6G_z$ is the projective connection
which satisfies
\begin{eqnarray}
&&
6G_t(t)dt^2=6G_z(z)dz^2+\{z,t\}dt^2
\nonumber
\end{eqnarray}
for another local coordinate $t$ around $P$,
where $\{z,t\}$ is the Schwarzian differential defined by
\begin{eqnarray}
&&
\{z,t\}={z^{'''}\over z^{'}}-{3\over2}\Big({z^{''}\over z^{'}}\Big)^2.
\nonumber
\end{eqnarray}

By Corollary \ref{expand2} and (\ref{remarkable}) we have

\begin{prop}\label{proj1}
\begin{eqnarray}
&&
G_z(z)=
{1\over N}
\sum_{i,j=1}^{Nm}
{q(k_i,k_j) dz(x)^2 \over (z(x)-\la_i)(z(x)-\la_j)}
-
\sum_{i,j=1}^g
{\partial^2 \log \theta[e_{\Lambda}] \over \partial z_i\partial z_j}
(0)v_i(x)v_j(x).
\nonumber
\end{eqnarray}
\end{prop}
\vskip2mm

As a corollary of this expression we have

\begin{cor}\label{projexp1}
Let $t=(z-\la_i)^{1/N}$ be the local coordinate around
the branch point $Q_i$.
Then the coefficient of $t^{N-2}dt$ in the Laurent expansion
of $G_z(z)$ in $t$ is
\begin{eqnarray}
&&
2N\sum_{j\neq i}{q(k_i,k_j)\over \la_i-\la_j}
%\nonumber
%\\
%&&
-{1\over (N-2)!}
\sum_{r,s=1}^g\sum_{\alpha=0}^{N-2}
\bic{N-2}{\alpha}
{\partial^2 \log \theta[e_{\Lambda}] \over \partial z_r \partial z_s}
(0)v_r^{(\alpha)}(Q_i)v_s^{(N-2-\alpha)}(Q_i),
\nonumber
\end{eqnarray}
where $v_r^{(\alpha)}(Q_i)$ is the coefficient of $dt$
in the expansion of $v_k(x)$ in $t$.
\end{cor}
\vskip1cm

\section{Another description of canonical differential}
\par

Let $P_l^{(l)}(z,w)$ be a polynomial
satisfying the conditions
\begin{description}
\item[1.] $P^{(l)}_l(z,w)=\sum_{j=0}^{lm}P_{lj}^{(l)}(w)(z-w)^j$ with
\begin{eqnarray}
&&
P_{l0}^{(l)}(w)=f(w),
\quad
P_{l1}^{(l)}(w)={l\over N}\sum_{i=1}^{Nm}{f(w)\over w-\la_i}.
\nonumber
\end{eqnarray}
\item[2.] $\hbox{deg}_wP_l^{(l)}(z,w)\leq (N-l)m$.
\end{description}
\vskip2mm

The following lemma can be easily proved.
\begin{lem}
The polynomial $P_l^{(l)}(z,w)$ saisfying the above conditions 1,2
exists for $l=1,\cdots,N-1$.
\end{lem}
\vskip2mm

We set
\begin{eqnarray}
\xi^{(0)}(x,y)&=&{dz(x)dz(y)\over(z(x)-z(y))^2},
\nonumber
\\
\xi^{(l)}(x,y)&=&
{P_l^{(l)}(z(x),z(y))dz(x)dz(y)\over s_l(x)s_{N-l}(y)(z(x)-z(y))^2}
\quad
\hbox{for } l=1,\cdots,N-1,
\nonumber
\\
\xi(x,y)&=&{1\over N}\sum_{l=0}^{N-1}\xi^{(l)}(x,y),
\nonumber
\end{eqnarray}
where $s_l(x)=s(x)^l$.
The condition 1,2 implies that $\xi^{(l)}(x,y)$
is regular on $C\times C$ except $\{z(x)=z(y)\}$.

\begin{prop}\label{expxi}

\begin{description}
\item[1.] $\xi(x,y)$ is holomorphic outside the diagonal set $\{x=y\}$.
\item[2.] For a non-branch point $P\in C$ we take $z$ 
as a local coordinate around $P$. Then the expansion in 
$z(x)$ at $z(y)$ is 
\begin{eqnarray}
&&
\xi(x,y)={dz(x)dz(y)\over(z(x)-z(y))^2}+O((z(x)-z(y))^0).
\nonumber
\end{eqnarray}
\end{description}
\end{prop}
\vskip2mm

\noindent
Proof. For $y\in C$ let $y^{(r)}=(z(y),\omega^rs(y))$.
Suppose that $y$ is not a branch point.
Then we can take $z$ as a local coordinate around $y^{(r)}$.
By calculation we have the expansion of $\xi^{(l)}(x,y)$
in $z(x)$ at $y^{(r)}$ as
\begin{eqnarray}
\xi^{(l)}(x,y)&=&
\omega^{-rl}\Big[
{1\over (z(x)-z(y))^2}
-{l^2\over 2N^2}\Big({d\over dz}\log f(z(y))\Big)^2
-{l\over 2N}{d^2\over dz^2}\log f(z(y))
\nonumber
\\
&&
+{P^{(l)}_{l,2}(z(y))\over f(z(y))}
+O\Big((z(x)-z(y))^1\Big)\Big]dz(x)dz(y).
\label{expxil}
\end{eqnarray}
By definition $\xi(x,y)$ is regular except $z(x)=z(y)$.
In order to prove the property 1 of the proposition
it is sufficient to prove that $\xi(x,y)$ has no singularity
at $x=y^{(r)}$ for $1\leq r\leq N-1$.
Note that if $y=Q_i$ for some $i$, then $z(x)=z(y)$ is equivalent to
$x=y=Q_i$. Hence by the expansion (\ref{expxil}), $\xi(x,y)$ is
regular at $x=y^{(r)}$. 
The property 2 is also obvious from (\ref{expxil}) above. $\Box$

\begin{cor}
$\omega(x,y)-\xi(x,y)$ is holomorphic on $C\times C$.
\end{cor}
\vskip2mm

\noindent
Proof. By Proposition \ref{expxi}, $\omega(x,y)-\xi(x,y)$
is regular except $\cup_{i=1}^{Nm}(Q_i,Q_i)$.
Hence
$\omega(x,y)-\xi(x,y)$ is regular everywhere on $C\times C$,
since $\omega(x,y)-\xi(x,y)$ is meromorphic on $C\times C$.
$\Box$
\vskip1mm

By this corollary there exists a set of polynomials
$P^{(l)}_k(z,w)$ such that
\begin{eqnarray}
&&
\omega(x,y)-\xi(x,y)=
\sum_{l=1}^{N-1}\sum_{k=1,k\neq l}^{N-1}
{P_k^{(l)}(z(x),z(y))dz(x)dz(y)\over s_k(x)s_{N-l}(y)},
\nonumber
\end{eqnarray}
where by changing the definition of $P^{(l)}_l(z,w)$
the $k=l$ term can be excluded.
The condition for the right hand side to be regular
at $z(x)=\infty$ and $z(y)=\infty$ is
\begin{eqnarray}
&&
\hbox{deg}_zP^{(l)}_k(z,w)\leq km-2,
\quad
\hbox{deg}_wP^{(l)}_k(z,w)\leq (N-l)m-2.
\nonumber
\end{eqnarray}
Hence we can write
\begin{eqnarray}
&&
P^{(l)}_k(z,w)=\sum_{j=0}^{km-2}P^{(l)}_{kj}(w)(z-w)^j,
\nonumber
\end{eqnarray}
for some polynomials $P^{(l)}_{kj}(w)$.
Now by the condition that the $A$ period of $\omega(x,y)$ is zero
we have

\begin{prop}
The following relation holds
\begin{eqnarray}
&&
\sum_{l=1}^{N-1}P^{(l)}_{l2}(\la_i)=
-f^\prime(\la_i){\partial\over\partial\la_i}\log\det A,
\nonumber
\end{eqnarray}
where $A$ is the $g\times g$ period matrix of non-normalized form:
\begin{eqnarray}
&&
A=(\int_{A_i}w^{(\alpha)}_\beta).
\nonumber
\end{eqnarray}
\end{prop}
\vskip2mm

\noindent
Proof. Let us take $t=(z-\la_i)^{1/N}$ as a local coordinate
around $Q_i$.
Then we have
\begin{eqnarray}
{dz(y)\over s_{N-l}(y)}&=&
{N\over\prod_{j\neq i}(\la_i-\la_j)^{(N-l)/N}}
t^{l-1}dt(1+O(t^N))
\quad 1\leq l\leq N-1,
\nonumber
\\
{dz(y)\over(z(x)-z(y))^2}
&=&
{N\over (z(x)-\la_i)^2}
t^{N-1}dt(1+O(t^N)).
\nonumber
\end{eqnarray}
Therefore if we set
\begin{eqnarray}
&&
\omega^{(l)}(x)=
{1\over N}
{P^{(l)}_l(z(x),\la_i)dz(x)\over s_l(x)(z(x)-\la_i)^2}
+\sum_{k=1,k\neq l}^{N-1}
{P^{(l)}_k(z(x),\la_i)dz(x)\over s_k(x)},
\nonumber
\end{eqnarray}
then the condition that
the coefficients of $dt,tdt,\cdots,t^{N-2}dt$ in the expansion of
$\int_{A_j}\omega(x,y)$ vanish
is equivalent to
\begin{eqnarray}
&&
\int_{A_j}\omega^{(l)}(x)=0
\quad
1\leq l\leq N-1.
\label{expomegal}
\end{eqnarray}
Noting that
\begin{eqnarray}
&&
P^{(l)}_l(z,\la_i)=
{l\over N}f^\prime(\la_i)(z-\la_i)
+\sum_{j=0}^{lm-2}P^{(l)}_{l,j+2}(\la_i)(z-\la_i)^{j+2},
\nonumber
\\
&&
{\partial\over\partial\la_i}{dz\over s_l}=
{l\over N}{dz\over s_l(z-\la_i)},
\nonumber
\end{eqnarray}
we see that (\ref{expomegal}) is equivalent to
\begin{eqnarray}
&&
{f^\prime(\la_i)\over N}{\partial\over\partial\la_i}
\int_{A_j}{dz\over s_l}+
{1\over N}\sum_{j=0}^{lm-2}P^{(l)}_{l,j+2}(\la_i)
\int_{A_j}{(z-\la_i)^jdz\over s_l}
\nonumber
\\
&&
+\sum_{k=1,k\neq l}\sum_{j=0}^{km-2}P^{(l)}_{k,j}(\la_i)
\int_{A_j}{(z-\la_i)^jdz\over s_k}=0.
\label{lineareq}
\end{eqnarray}
We consider (\ref{lineareq}) as a linear equation for
the $g$ variables $\{P^{(l)}_{k,r}(\la_i)\}$.
Solving (\ref{lineareq}) in $P^{(l)}_{l,2}$ 
by the Cramer's formula and summing up in $l$ we
have the statement of the proposition. $\Box$
\vskip2mm

The idea of deriving equations of the form (\ref{lineareq})
is due to Bershadsky-Radul\cite{BR1}.
By calculations we have
\begin{cor}\label{projexp2}
The coefficient of $t^{N-2}dt$ in the expansion of $G_z(z)$ in
$t=(z-\la_i)^{1/N}$ is
\begin{eqnarray}
&&
-\mu N\sum_{j=1,j\neq i}^{Nm}{1\over \la_i-\la_j}-
N{\partial\over\partial\la_i}\log\det A,
\nonumber
\end{eqnarray}
where
\begin{eqnarray}
&&
\mu={(N-1)(2N-1)\over 6N}.
\nonumber
\end{eqnarray}
\end{cor}
\vskip1cm

\section{Variational formula of period matrix}
\par
Let us consider the equation
\begin{eqnarray}
&&
s_t^N=(z-\la_i-t)\prod_{j=1,j\neq i}^{Nm}(z-\la_j)
\nonumber
\end{eqnarray}
which is a one parameter deformation of the curve $C$
by a small parameter $t$. 
We denote the corresponding compact Riemann surface by $C_t$.
The notation $s_t$ is different from $s_l=s^l$ in the previous section.
We hope that this does not cause any confusion.
Let $\bar{\pi}$ be the projection $\bar{\pi}:C\ra {\bf P}^1$
which maps $(z,s)$ to $z$.
We can take a canonical dissection $\{A_i(t),B_j(t)\}$ of $C_t$
such that $\bar{\pi}(A_i(t))$, $\bar{\pi}(B_j(t))$ do not depend
on $t$ for $\vert t\vert$ being sufficiently small.
The integration of a holomorphic $1$-form on $C_t$
along $A_i(t),B_j(t)$ can be considered as the integration
of a multi-valued holomorphic $1$-form on 
${\bf P}^1-\{\la_1,\cdots,\la_{Nm}\}$ along
$\bar{\pi}(A_i(t))$, $\bar{\pi}(B_j(t))$.
Hence we can think of the integration cycles $A_i(t),B_j(t)$
as if they are independent of $t$.
Therefore we simply write $A_i,B_j$ instead of
$A_i(t),B_j(t)$ in the calculations in this section.
Let $\{v_j(x,t)\}$ be the basis of 
normailzed holomorphic 1-forms on $C_t$ with respect to 
$\{A_i(t),B_j(t)\}$.
We denote by $\tau(t)=(\tau_{kr}(t))$ the period matrix
\begin{eqnarray}
&&
\tau_{kr}(t)=\int_{B_k}v_r(x,t).
\nonumber
\end{eqnarray}
We set
\begin{eqnarray}
&&
w^{(\alpha)}_{\beta t}={z^{\beta-1}dz\over s_t^\alpha}.
\nonumber
\end{eqnarray}
We also use our previous notation $v_j(x)=v_j(x,0)$, $\tau_{kr}=\tau_{kr}(0)$,
$s=s_0$, $w^{(\alpha)}_\beta=w^{(\alpha)}_{\beta0}$.
Our aim in this section is to prove

\begin{theo}\label{variation}
\begin{eqnarray}
&&
{d\tau_{jk}\over dt}(0)=
{1\over N(N-2)!}\sum_{\alpha=0}^{N-2}
\bic{N-2}{\alpha}v_j^{(\alpha)}(Q_i)v_k^{(N-2-\alpha)}(Q_i).
\nonumber
\end{eqnarray}
\end{theo}
\vskip2mm

We define the connection matrix $\sigma$ and $c$ by
\begin{eqnarray}
&&
v_j(x)=\sum_{\alpha,\beta}\sigma_{j(\alpha\beta)}w^{(\alpha)}_\beta(x),
\quad
w^{(\alpha)}_\beta(x)=\sum_{j}c_{(\alpha\beta)j}v_j(x).
\label{connection}
\end{eqnarray}

Let $P\in C$ and $u$ be a local coordinate around $P$.
Let $\omega(P;n)$ be 
the abelian differential of the second kind
satisfying the following conditions.
\vskip2mm

\begin{description}
\item[1.] $\omega(P;n)$ is holomorphic except the point $P\in C$
where it has a pole of order $n\geq2$. 
At $P$ we have the expansion of the form
\begin{eqnarray}
&&
\omega(P;n)=-{n-1\over u^n}du(1+O(u^n)).
\nonumber
\end{eqnarray}
\item[2.] $\omega(P;n)$ has zero $A$ periods :
\begin{eqnarray}
&&
\int_{A_j}\omega(P;n)=0
\quad\hbox{ for any $j$.}
\nonumber
\end{eqnarray}
\end{description}
\vskip2mm

The differential $\omega(P;n)$ 
depends on the choice of the local coordinate $u$.
In our case $P=Q_i$ we always take $u=(z-\la_i)^{1/N}$ as a local
coordinate around $Q_i$. In this sense $\omega(P;n)$ is uniquely determined.
It is known that the following relation holds
\begin{eqnarray}
&&
\int_{B_j}\omega(P;n)=-{1\over(n-2)!}v_j^{(n-2)}(P),
\label{nicerels}
\end{eqnarray}
where $v_j^{(n-2)}(P)$ is the coefficient of $u^{n-2}du$ in the expansion
of $v_j(x)$ in $u$.

\begin{lem}
If we expand $v_j(x,t)$ as
\begin{eqnarray}
&&
v_j(x,t)=v_j(x)+v_{j1}(x)t+\cdots,
\label{fexpvt}
\end{eqnarray}
then we have
\begin{eqnarray}
&&
v_{j1}(x)=
-\sum_{\alpha,\beta}
{\sigma_{j(\alpha\beta)}\la_i^{\beta-1}\over
\prod_{j\neq i}(\la_i-\la_j)^{\alpha/N}}
\omega(Q_i;\alpha+1).
\nonumber
\end{eqnarray}
\end{lem}
\vskip2mm

\noindent
Proof. We have the expansion
\begin{eqnarray}
&&
w^{(\alpha)}_{\beta t}(x)=
w^{(\alpha)}_\beta(x)+{\alpha z^{\beta-1}dz\over N(z-\la_i)s^\alpha}t
+O(t^2),
\label{expwt}
\end{eqnarray}
and the relation
\begin{eqnarray}
w^{(\alpha)}_{\beta t}(x)=
\sum_{j=1}^g\int_{A_j}w^{(\alpha)}_{\beta t}\cdot v_j(x,t).
\label{relwv}
\end{eqnarray}
Substituting the expansions (\ref{fexpvt}) and (\ref{expwt})
into the equation (\ref{relwv}) and comparing the coefficient
of $t$ we have
\begin{eqnarray}
\sum_{j=1}^gv_{j1}(x)\int_{A_j}w^{(\alpha)}_\beta&=&
{\alpha\over N}\eta^{(\alpha)}_\beta(x),
\nonumber
\\
\eta^{(\alpha)}_\beta(x)&=&
{z^{\beta-1}dz\over(z-\la_i)s^\alpha}
-\sum_{j=1}^gv_j(x)\int_{A_j}{z^{\beta-1}dz\over(z-\la_i)s^\alpha}.
\label{exprels}
\end{eqnarray}
Then $\eta^{(\alpha)}_\beta(x)$ has the following properties:
\vskip2mm
\begin{description}
\item[1.] $\int_{A_k}\eta^{(\alpha)}_\beta(x)=0$ for any $k=1,\cdots,g$.
\item[2.] Taking $u=(z-\la_i)^{1/N}$ as a local coordinate around $Q_i$
we have the expansion
\begin{eqnarray}
&&
\eta^{(\alpha)}_\beta(x)=
{N\la_i^{\beta-1}\over f^\prime(\la_i)^{\alpha/N}}
{du\over u^{\alpha+1}}+O(1).
\nonumber
\end{eqnarray}
\end{description}
\vskip2mm

Hence we have
\begin{eqnarray}
&&
\eta^{(\alpha)}_\beta(x)
=-{N\la_i^{\beta-1}\over \alpha f^\prime(\la_i)^{\alpha/N}}
\omega(Q_i;\alpha+1).
\label{etacanonical}
\end{eqnarray}
Since
\begin{eqnarray}
&&
\int_{A_j}w^{(\alpha)}_\beta=c_{(\alpha\beta)j},
\nonumber
\end{eqnarray}
and $\sigma$ is the inverse matrix of $c$, we have the
desired result from (\ref{exprels}) and (\ref{etacanonical}).
$\Box$
\vskip2mm

Now comparing the coefficient of $u^{N-1-\alpha}du$
of the both hand sides of the first equation of 
(\ref{connection}) we have
\begin{eqnarray}
&&
\sum_{\beta=1}^{\alpha m-1}\sigma_{j(\alpha\beta)}\la_i^{\beta-1}
={f^\prime(\la_i)^{\alpha/N}\over
N(N-1-\alpha)!}
v_j^{(N-1-\alpha)}(Q_i),
\nonumber
\end{eqnarray}
for $1\leq \alpha\leq N-1$ and thus
\begin{eqnarray}
v_{j1}(x)=
-{1\over N}\sum_{\alpha=1}^{N-1}
{1\over (N-1-\alpha)!}
v_j^{(N-1-\alpha)}(Q_i)\omega(Q_i;\alpha+1).
\nonumber
\end{eqnarray}
Integrating both hand sides of this equation along the cycle
$B_k$ and using the relation (\ref{nicerels}) we obtain
\begin{eqnarray}
&&
\int_{B_k}v_{j1}(x)=
{1\over N(N-2)!}
\sum_{\alpha=0}^{N-2}
\bic{N-2}{\alpha}v_j^{(N-2-\alpha)}(Q_i)v_k^{(\alpha)}(Q_i).
\nonumber
\end{eqnarray}
$\Box$
\vskip1cm

\section{Thomae formula}
\par
Now let us prove the generalized Thomae formula.

\begin{theo}\label{Thomae}
For an ordered partition $\La=(\La_0,\cdots,\La_{N-1})$ we have
\begin{eqnarray}
&&
\theta[e_\La](0)^{2N}
=C_{\La}(\det A)^{N}\prod_{i<j}(\la_i-\la_j)^{2Nq(k_i,k_j)+N\mu},
\nonumber
\end{eqnarray}
where $k_i=j$ for $i\in\La_j$,
\begin{eqnarray}
&&
q(i,j)=\sum_{l\in{\cal L}}q_l(i)q_l(j),
\quad
\mu={(N-1)(2N-1)\over 6N},
\nonumber
\end{eqnarray}
and $q_l(i)$, ${\cal L}$ are given by (\ref{qli}), (\ref{calL})
in section 3.
The complex number
$C_{\La}$ does not depend on $\la_i$'s.
They satisfy
$C_{\La}^{2N}=C_{\La^\prime}^{2N}$ for any $\La$, $\La^\prime$.
\end{theo}
\vskip3mm

Since the familly $\{C_t\}$ is locally topologically trivial,
we can take a canonical dissection $\{A_i(t), B_j(t)\}$ of $C_t$
such that $A_i(t)$, $B_j(t)$ are continuous in $t$.
We assume that $A_i(t),B_j(t)$ does not go through any
branch point $Q_i(t)$.
We can also define the base points $P_0(t)$ of $C_t$ and
$z_0(t)$ of $\tilde{C}_t$ lying over $P_0(t)$ 
so that they vary continuously in $t$.
We identify $Q_i(t)$ with the corresponding point in the
fundamental domain on $\tilde{C}_t$ which contains the base point
$z_0(t)$. 
Let $k^{P_0}(t)$ be the vector in ${\bf C}^g$ whose $j$-th
component is defined by
\begin{eqnarray}
&&
k^{P_0}(t)_j=
{2\pi i-\tau_{jj}(t)\over2}
+{1\over 2\pi i}\sum_{i\neq j}
\int_{A_i(t)}v_i(x,t)\int_{z_0(t)}^xv_j(x,t).
\nonumber
\end{eqnarray}
It is known (see \cite{F}(p8) for example) that
\begin{eqnarray}
&&
\Delta-(g-1)P_0(t)=k^{P_0}(t)
\nonumber
\end{eqnarray}
in $J(C)$.
Let us define $e_\La(t)$ as an element of ${\bf C}^g$ by
\begin{eqnarray}
&&
e_\La(t)=
\sum_{j\in\La_1}\int_{z_0(t)}^{Q_j(t)}v(x,t)+\cdots
+(N-1)\sum_{j\in\La_{N-1}}\int_{z_0(t)}^{Q_j(t)}v(x,t)
-k^{P_0}(t).
\nonumber
\end{eqnarray}
Then $e_\La(t)$ is continuous in $t$.
Note that the linear isomorphism ${\bf C}^g\simeq {\bf R}^{2g}$
sending $e$ to its characteristics with respect $\tau(t)$ is
analytic in $t$. Therefore if we write
\begin{eqnarray}
&&
e_\La(t)=\cht{\delta(t)}{\ep(t)},
\nonumber
\end{eqnarray}
then $\ep(t)$ and $\delta(t)$ are continuous in $t$.
Since $\ep(t)$ and $\delta(t)$ are in $1/2N{\bf Z}^{g}$,
they are constant in $t$.
Therefore we simply write $\ep,\delta$ instead of $\ep(t),\delta(t)$.
We denote by $\theta_t[e](z)$ the theta function
associated with the canonical basis $\{A_i(t),B_j(t)\}$
of $C_t$. We set $\theta[e](z)=\theta_0[e](z)$.
Then the function $\theta_t[e_\La(t)](0)$ depends on $t$ 
only through the period matrix $\tau_{kr}(t)$ since
\begin{eqnarray}
&&
\theta_t[e_\La(t)](0)=
\sum_{m\in\bzg}\exp({1\over2}(m+\delta)\tau(t)(m+\delta)^t
+2\pi i\ep(m+\delta)^t).
\nonumber
\end{eqnarray}
Using the heat equations
\begin{eqnarray}
&&
{\partial^2\theta[e_\La]\over \partial z_k\partial z_r}(z)
={\partial\theta[e_\La]\over\partial\tau_{kr}}(z)
\quad(k\neq r),
\quad
{\partial^2\theta[e_\La]\over \partial z_k^2}(z)
=2{\partial\theta[e_\La]\over\partial\tau_{kk}}(z),
\nonumber
\end{eqnarray}
and Lemma \ref{vanish} we have
\begin{eqnarray}
{\partial\over\partial\la_i}\log\theta[e_\La](0)&=&
{d\over dt}\log\theta_t[e_\La(t)](0)\Big\vert_{t=0}
\nonumber
\\
&=&\sum_{1\leq k\leq r\leq g}
{\partial \log\theta[e_\La]\over \partial\tau_{kr}}(0)
{d\tau_{kr}\over dt}(0)
\nonumber
\\
&=&{1\over2}
\sum_{k,r=1}^g
{1\over \theta[e_\La](0)}
{\partial^2 \theta[e_\La]\over \partial z_k\partial z_r}(0)
{d\tau_{kr}\over dt}(0),
\nonumber
\\
&=&
{1\over2}\sum_{k,r=1}^g
{\partial^2 \log\theta[e_\La]\over \partial z_k\partial z_r}(0)
{d\tau_{kr}\over dt}(0).
\label{logderiv}
\end{eqnarray}

On the other hand by Corollary \ref{projexp1}, \ref{projexp2} and
Theorem \ref{variation} we have
\begin{eqnarray}
&&-\mu N\sum_{j\neq i}{1\over\la_i-\la_j}
-N{\partial\over\partial\la_i}\log\det A
\nonumber
\\
&&=
2N\sum_{j\neq i}{q(k_i,k_j)\over\la_i-\la_j}
-N\sum_{k,r=1}^g
{\partial^2 \log\theta[e_\La]\over \partial z_k\partial z_r}(0)
{d\tau_{kr}\over dt}(0).
\label{preeq}
\end{eqnarray}

Substituting (\ref{logderiv}) into (\ref{preeq}) we have

\begin{eqnarray}
&&{\partial\over\partial\la_i}\log\theta[e_\La](0)
=
{1\over2}{\partial\over\partial\la_i}\log\det A
+{\mu\over2}\sum_{j\neq i}{1\over\la_i-\la_j}
+\sum_{j\neq i}{q(k_i,k_j)\over\la_i-\la_j}.
\nonumber
\end{eqnarray}
Hence we have proved the first part of Theorem \ref{Thomae}.
\vskip1cm

\section{Property of the constant $C_\La$}
\par
Our aim in this section is to prove 
the remaining part of Theorem \ref{Thomae}, that is,
for any ordered partitions $\La$ and $\La^\prime$
\begin{eqnarray}
&&
C_{\La}^{2N}=C_{\La^\prime}^{2N}.
\label{remaining}
\end{eqnarray}
As in the previous section we identify branch points $Q_i$
with the corresponding points in the fundamental domain in $\tilde{C}$.

The key for the proof is the formula of Fay ( \cite{F}, p30, Cor.2.17 ):
\begin{eqnarray}
&&
\theta(\sum_{k=1}^gx_k-p-\Delta)
=c{\det(v_i(x_j))\over\prod_{i<j}E(x_i,x_j)}
{\sigma(p)\over \prod_{k=1}^g\sigma(x_k)}
\prod_{k=1}^gE(x_k,p),
\nonumber
\end{eqnarray}
for any $p,x_1,\ldots,x_g\in \tilde{C}$, where $c$ is independent on
$p,x_1,\ldots,x_g$ and
\begin{eqnarray}
&&
\sigma(p)=
\exp(-{1\over2\pi i}\sum_{j=1}^g\int_{A_j}v_j(y)\log E(y,p)).
\nonumber
\end{eqnarray}
Taking ratios for $p=a,b$ equations we have
\begin{eqnarray}
&&
{\theta(\sum_{k=1}^gx_k-a-\Delta)\over \theta(\sum_{k=1}^gx_k-b-\Delta)}
=
{\sigma(a)\over\sigma(b)}
\prod_{k=1}^g
{E(x_k,a)\over E(x_k,b)}.
\nonumber
\end{eqnarray}
Set $a=Q_i$, $b=Q_j$ $(i\neq j)$ and taking $N$-th power of the
both hand sides we obtain
\begin{eqnarray}
&&
{\theta(\sum_{k=1}^gx_k-Q_i-\Delta)^N
\over \theta(\sum_{k=1}^gx_k-Q_j-\Delta)^N}
=
\Big({\sigma(Q_i)\over\sigma(Q_j)}\Big)^N
\prod_{k=1}^g
{E(x_k,Q_i)^N\over E(x_k,Q_j)^N}.
\nonumber
\end{eqnarray}
\noindent
Since $NQ_i$ and $NQ_j$ are linearly equivalent,
there exists $\la(i,j),\kappa(i,j)\in\bzg$ such that
\begin{eqnarray}
&&
N\int_{z_0}^{Q_i}v(x)-N\int_{z_0}^{Q_j}v(x)
=N\int_{Q_i}^{Q_j}v(x)=2\pi i\la(i,j)+\kappa(i,j)\tau,
\nonumber
\\
&&\kappa(j,i)=-\kappa(i,j),\quad
\la(j,i)=-\la(i,j).
\nonumber
\end{eqnarray}
If we set 
\begin{eqnarray}
&&
f(x)=
{E(x,Q_i)^N\over E(x,Q_j)^N},
\nonumber
\end{eqnarray}
then we have
\begin{eqnarray}
f(x+A_k)&=&f(x)
\nonumber
\\
f(x+B_k)&=&
\exp(-\sum_{l}\tau_{lk}\kappa(i,j)_l)f(x).
\nonumber
\end{eqnarray}
Hence the function
\begin{eqnarray}
&&
\exp(\sum_{l}\int_{z_0}^xv_l(x)\kappa(i,j)_l)f(x)
\nonumber
\end{eqnarray}
can be considered as a single valued function on $C$.
Its only zeros are of $N$-th order at $Q_i$
and only poles are of $N$-th order at $Q_j$.
Therefore there exists a constant $c_{ij}$ such that
\begin{eqnarray}
&&
\exp(\int_{z_0}^xv(x)\kappa(i,j)^t)
{E(x,Q_i)^N\over E(x,Q_j)^N}
=c_{ij}{z(x)-\la_i\over z(x)-\la_j}.
\nonumber
\end{eqnarray}
By the property of $\kappa(i,j)$, $c_{ij}$ satisfies
\begin{eqnarray}
&&
c_{ij}=c_{ji}^{-1},
\quad
c_{ii}=1.
\nonumber
\end{eqnarray}
If we set
\begin{eqnarray}
&&
r_{ij}=c_{ij}\Big({\sigma(Q_i)\over\sigma(Q_j)}\Big)^N,
\nonumber
\\
&&
w(x\vert i,j)=\int_{z_0}^xv(x)\kappa(i,j)^t,
\quad
w(\sum_kx_k\vert i,j)=\sum_k w(x_k \vert i,j),
\nonumber
\end{eqnarray}
we have
\begin{eqnarray}
&&
\exp(w(\sum_{k=1}^gx_k\vert i,j))
{\theta(\sum_{k=1}^gx_k-Q_i-\Delta)^N
\over \theta(\sum_{k=1}^gx_k-Q_j-\Delta)^N}
=r_{ij}\prod_{k=1}^g
{z(x)-\la_i\over z(x)-\la_j},
\label{ratio1}
\\
&&
r_{ij}=r_{ji}^{-1},
\quad
r_{ii}=1.
\nonumber
\end{eqnarray}

Now let us take an ordered partition
$\La=\La^{(1)}=(\La^{(1)}_{0},\ldots,\La^{(1)}_{N-1})$ with
\begin{eqnarray}
&&
\La^{(1)}_l=\{i^l_1,\ldots,i^l_m\},
\quad
0\leq l\leq N-1.
\nonumber
\end{eqnarray}
Let us define $\La^{(2)}=(\La^{(2)}_{0},\ldots,\La^{(2)}_{N-1})$
by
\begin{eqnarray}
&&
\La^{(2)}_0=\{i^0_1,\ldots,i^0_{m-1},i^{N-1}_m\},
\quad
\La^{(2)}_{N-1}=\{i^{N-1}_1,\ldots,i^{N-1}_{m-1},i^{0}_m\},
\nonumber
\end{eqnarray}
and $\La^{(2)}_l=\La^{(2)}_l(1\leq l\leq N-1)$.

If we consider $Q_i$ as $\int_{z_0}^{Q_i}v$,
we have the vectors in $\bcg$:
\begin{eqnarray}
e_{\La^{(1)}}&=&\La^{(1)}_1+2\La^{(1)}_2+\cdots+
(N-1)(\La^{(1)}_{N-1}\backslash\{i^{N-1}_m\})-Q_{i^{N-1}_m}-
k^{P_0},
\nonumber
\\
e_{\La^{(2)}}&=&\La^{(2)}_1+2\La^{(2)}_2+\cdots+
(N-1)(\La^{(2)}_{N-1}\backslash\{i^{0}_m\})-Q_{i^{0}_m}-k^{P_0},
\nonumber
\\
&=&\La^{(1)}_1+2\La^{(1)}_2+\cdots+
(N-1)(\La^{(1)}_{N-1}\backslash\{i^{N-1}_m\})-Q_{i^{0}_m}-k^{P_0}.
\nonumber
\end{eqnarray}
Putting $i=i^{N-1}_m$, $j=i^0_m$ and
\begin{eqnarray}
&&
(x_1,\ldots,x_g)=(\La^{(1)}_1,2\La^{(1)}_2,\ldots,
(N-1)(\La^{(1)}_{N-1}\backslash\{i^{N-1}_m\}))
\nonumber
\end{eqnarray}
in (\ref{ratio1}) we have
\begin{eqnarray}
&&
\exp(U(\La^{(1)}\vert i^{N-1}_m,i^0_m))
{\theta(e_{\La^{(1)}})^N\over \theta(e_{\La^{(2)}})^N}
=r_{i^{N-1}_m,i^0_m}
\prod_{r=1}^{N-1}
{\prod_{s=1}^{m}}^\prime
\Big(
{\la_{i^r_s}-\la_{i^{N-1}_m}\over \la_{i^r_s}-\la_{i^{0}_m}}
\Big)^r,
\label{ratio2}
\\
&&
U(\La^{(1)}\vert i^{N-1}_m,i^0_m)=
w(\La^{(1)}_1+2\La^{(1)}_2+\cdots+
(N-1)(\La^{(1)}_{N-1}\backslash\{i^{N-1}_m\})\vert
i^{N-1}_m,i^0_m).
\nonumber
\end{eqnarray}
Here $\prod'$ means the product for $(r,s)\neq(N-1,m)$.

Let us define the elements of $\bcg$ by
\begin{eqnarray}
-\bar{e}_{\La^{(1)}}&=&
\La^{(1)}_{N-2}+2\La^{(1)}_{N-3}+\cdots+
(N-1)(\La^{(1)}_{0}\backslash\{i^{0}_m\})
-Q_{i^{0}_m}-k^{P_0},
\nonumber
\\
-\bar{e}_{\La^{(2)}}&=&
\La^{(2)}_{N-2}+2\La^{(2)}_{N-3}+\cdots+
(N-1)(\La^{(2)}_{0}\backslash\{i^{N-1}_m\})
-Q_{i^{N-1}_m}-k^{P_0}
\nonumber
\\
&=&
\La^{(1)}_{N-2}+2\La^{(1)}_{N-3}+\cdots+
(N-1)(\La^{(1)}_{0}\backslash\{i^{0}_m\})
-Q_{i^{N-1}_m}-k^{P_0},
\nonumber
\end{eqnarray}
where again $Q_i$ denotes $\int_{z_0}^{Q_i}v$.
Then if we set
$i=i^{0}_m$, $j=i^{N-1}_m$ and
\begin{eqnarray}
&&
(x_1,\ldots,x_g)=(\La^{(1)}_{N-2},2\La^{(1)}_{N-3},\ldots,
(N-1)(\La^{(1)}_{0}\backslash\{i^{0}_m\}))
\nonumber
\end{eqnarray}
in (\ref{ratio1}) we have
\begin{eqnarray}
&&
\exp(U'(\La^{(1)}\vert i^{0}_m,i^{N-1}_m))
{\theta(\bar{e}_{\La^{(1)}})^N\over \theta(\bar{e}_{\La^{(2)}})^N}
=r_{i^0_m,i^{N-1}_m}
\prod_{r=1}^{N-1}
{\prod_{s=1}^{m}}^\prime
\Big(
{\la_{i^{N-1-r}_s}-\la_{i^{0}_m}\over \la_{i^{N-1-r}_s}-\la_{i^{N-1}_m}}
\Big)^r,
\label{ratio3}
\\
&&
U'(\La^{(1)}\vert i^{0}_m,i^{N-1}_m)=
w(\La^{(1)}_{N-2}+2\La^{(1)}_{N-3}+\cdots+
(N-1)(\La^{(1)}_{0}\backslash\{i^{0}_m\})\vert
i^{0}_m,i^{N-1}_m).
\nonumber
\end{eqnarray}

Here we have used the property that $\theta(z)$ is an even function of $z$.
Multiplying (\ref{ratio2}) and (\ref{ratio3}) we have

\begin{eqnarray}
&&
\exp\big(w(-\sum_{l=0}^{N-1}(N-1-2l)
\tilde{\La}^{(1)}_l\vert i^{N-1}_m,i^0_m)\big)
\Big({
\theta(e_{\La^{(1)}})\theta(\bar{e}_{\La^{(1)}})
\over 
\theta(e_{\La^{(2)}})\theta(\bar{e}_{\La^{(2)}})
}\Big)^N
\nonumber
\\
&&=\prod_{r=1}^{N-1}{\prod_{s=1}^{m}}^\prime
{
(\la_{i^r_s}-\la_{i^{N-1}_m})^r(\la_{i^{N-1-r}_s}-\la_{i^{0}_m})^r
\over 
(\la_{i^r_s}-\la_{i^{0}_m})^r(\la_{i^{N-1-r}_s}-\la_{i^{N-1}_m})^r
},
\label{ratio4}
\end{eqnarray}
where we set
\begin{eqnarray}
&&
\tilde{\La}^{(1)}_0=\La^{(1)}_0\backslash\{i^0_m\},
\quad
\tilde{\La}^{(1)}_{N-1}=\La^{(1)}_{N-1}\backslash\{i^{N-1}_m\},
\quad
\tilde{\La}^{(1)}_l=\La^{(1)}_l\quad(l\neq 0,N-1).
\nonumber
\end{eqnarray}
Since $\bar{e}_{\La^{(k)}}=e_{\La^{(k)}}$, $k=1,2$ in $J(C)$
we can set
\begin{eqnarray}
&&
e_{\La^{(k)}}=\ch{\delta^{(k)}}{\ep^{(k)}}=
2\pi i\ep^{(k)}+\delta^{(k)}\tau,
\quad \delta^{(k)},\ep^{(k)}\in{1\over 2N}\bzg,
\nonumber
\\
&&
\bar{e}_{\La^{(k)}}=e_{\La^{(k)}}+2\pi im^{(k)}+n^{(k)}\tau,
\quad m^{(k)},n^{(k)}\in\bzg.
\nonumber
\end{eqnarray}
Substituting these equations into (\ref{ratio4}), 
taking $2N$-th power of both hand sides and using the transformation
property of theta functions, we get
\begin{eqnarray}
&&
\Big({
\theta[e_{\La^{(1)}}](0)
\over
\theta[e_{\La^{(2)}}](0)
}\Big)^{4N^2}
=B
\prod_{r=1}^{N-1}{\prod_{s=1}^{m}}^\prime
{
(\la_{i^r_s}-\la_{i^{N-1}_m})^{2Nr}
(\la_{i^{N-1-r}_s}-\la_{i^{0}_m})^{2Nr}
\over 
(\la_{i^r_s}-\la_{i^{0}_m})^{2Nr}
(\la_{i^{N-1-r}_s}-\la_{i^{N-1}_m})^{2Nr}
},\label{ratio5}
\\
&&
B=
\exp\big(
-2n^\prime\tau\kappa^t-
N^2\sum_{k=1}^2(-1)^k(
n^{(k)}\tau n^{(k)t}
+2\delta^{(k)}\tau n^{(k)t}
+2\delta^{(k)}\tau \delta^{(k)t})
\big),
\nonumber
\end{eqnarray}
where we set
\begin{eqnarray}
&&
Nw(-\sum_{l=0}^{N-1}(N-1-2l)\tilde{\La}^{(1)}_l\vert i^{N-1}_m,i^0_m)
=(2\pi im^\prime+n^\prime\tau)\kappa,
\nonumber
\\
&&
\kappa=\kappa(i^{N-1}_m,i^0_m).
\nonumber
\end{eqnarray}

We can simplify the right hand side of (\ref{ratio5})
so that there are no common divisor in the numerator and the denominator.
The result is
\begin{eqnarray}
\Big({
\theta[e_{\La^{(1)}}](0)
\over
\theta[e_{\La^{(2)}}](0)
}\Big)^{4N^2}=
&&
B
\prod_{r=1}^{N-2}{\prod_{s=1}^{m}}
\Big({
\la_{i^r_s}-\la_{i^{0}_m}
\over
\la_{i^r_s}-\la_{i^{N-1}_m}
}\Big)^{2N(N-1-2r)}
\nonumber
\\
&&
\prod_{s=1}^{m-1}
\Big({
(\la_{i^{N-1}_s}-\la_{i^{N-1}_m})
(\la_{i^0_s}-\la_{i^{0}_m})
\over
(\la_{i^0_s}-\la_{i^{N-1}_m})
(\la_{i^{N-1}_s}-\la_{i^{0}_m})
}\Big)^{2N(N-1)}.
\label{ratio6}
\end{eqnarray}

Let us compare this equation with those obtained from
the proved part of Theorem \ref{Thomae}.
Let $\{k_i\}$, $\{k_i^\prime\}$ correspond to
$\La^{(1)}$, $\La^{(2)}$ respectively
as in (\ref{weight}). Then by the proved part of the
Thomae formula we have
\begin{eqnarray}
\Big({
\theta[e_{\La^{(1)}}](0)
\over
\theta[e_{\La^{(2)}}](0)
}\Big)^{4N^2}
=
C\prod_{i<j}(\la_i-\la_j)^{4N^2(q(k_i,k_j)-q(k^\prime_i,k^\prime_j))},
\label{conseq}
\end{eqnarray}
where $C=\big(C_{\La^{(1)}}/C_{\La^{(2)}}\big)^{4N^2}$.
Note that $4N^2q(k_i,k_j)$ and $4N^2q(k^\prime_i,k^\prime_j)$
are even.
Then we have
\begin{eqnarray}
&&
\hbox{LHS of (\ref{conseq})}
\nonumber
\\
&&
=
C 
\prod_{r=1}^{N-1}{\prod_{s=1}^{m}}^\prime
\Big({
\la_{i^r_s}-\la_{i^{N-1}_m}
\over 
\la_{i^r_s}-\la_{i^{0}_m}
}\Big)^{4N^2q(r,N-1)}
\prod_{r=0}^{N-2}{\prod_{s=1}^{m}}^{''}
\Big({
\la_{i^r_s}-\la_{i^{0}_m}
\over 
\la_{i^r_s}-\la_{i^{N-1}_m}
}\Big)^{4N^2q(r,0)}
\nonumber
\\
&&
\prod_{s=1}^{m-1}
\Big({
(\la_{i^{N-1}_s}-\la_{i^{0}_m})
(\la_{i^{0}_s}-\la_{i^{N-1}_m})
\over
(\la_{i^{N-1}_s}-\la_{i^{N-1}_m})
(\la_{i^{0}_s}-\la_{i^{0}_m})
}
\Big)^{4N^2q(0,1)}
\nonumber
\\
&&
=C
\prod_{r=1}^{N-2}\prod_{s=1}^m
\Big({
\la_{i^r_s}-\la_{i^{0}_m}
\over 
\la_{i^r_s}-\la_{i^{N-1}_m}
}\Big)^{4N^2(q(r,0)-q(r,N-1))}
\nonumber
\\
&&
\prod_{s=1}^{m-1}
\Big({
(\la_{i^{N-1}_s}-\la_{i^{N-1}_m})
(\la_{i^{0}_s}-\la_{i^{0}_m})
\over
(\la_{i^{N-1}_s}-\la_{i^{0}_m})
(\la_{i^{0}_s}-\la_{i^{N-1}_m})
}
\Big)^{4N^2(q(0,0)-q(0,1))},
\label{ratiothoma}
\end{eqnarray}
where  $\prod''$ means the product for $(r,s)\neq(0,m)$.
Let us calculate $q_l(r,t)$.
By the definition of $q(i,j)$ and Lemma \ref{exponent},
$q(i,j)$ depends only on $\vert i-j\vert$ mod $N$.
Hence using $q_l(0)=N/l$ we have
\begin{eqnarray}
&&
q(r,t)=q(0,t-r)
=\sum_{l\in{\cal L}}q_l(0)q_l(t-r)
={1\over N}\sum_{l\in{\cal L}}lq_l(t-r).
\nonumber
\end{eqnarray}
The following lemma is obtained by a direct calculation.

\begin{lem}
\begin{eqnarray}
\sum_{l\in{\cal L}}lq_l(r)={N^2-1\over12}-{1\over2}r(N-r).
\end{eqnarray}
\end{lem}
\vskip2mm

Thus we have
\begin{eqnarray}
&&
q(r,t)=
{1\over N}
\Big({N^2-1\over 12}-{1\over2}(t-r)(N-t+r)\Big).
\nonumber
\end{eqnarray}
In particular
\begin{eqnarray}
4N^2\big(q(0,0)-q(0,1)\big)&=&2N(N-1),
\nonumber
\\
4N^2\big(q(r,0)-q(r,N-1)\big)&=&2N(N-1-2r).
\nonumber
\end{eqnarray}
Comparing (\ref{ratio6}) and (\ref{ratiothoma}) we have
\begin{eqnarray}
\Big({C_{\La^{(1)}}\over C_{\La^{(2)}}}\Big)^{2N}
=B.
\nonumber
\end{eqnarray}
Let us write
\begin{eqnarray}
&&
B=\exp(\sum_{k\leq r}b_{kr}\tau_{kr}).
\nonumber
\end{eqnarray}
Since $C_{\La^{(1)}}$ and $C_{\La^{(2)}}$ do not depend on
$\tau$ 
\begin{eqnarray}
&&
{\partial\over\partial \tau_{kr}}B=b_{kr}B=0,
\quad\hbox{for any $k\leq r$}.
\nonumber
\end{eqnarray}
Hence $B=1$.
Since any two ordered partitions are transformed to each other
by successive exchange of elements of $\La_i$ and 
$\La_{i+1}$,$i=0,\ldots,N-1$,
the equation (\ref{remaining}) are proved.
\vskip1cm

\section{Examples}
\par
In this section we shall give examples of
Thomae formula for small $N$'s.
Recall that 
\begin{eqnarray}
&&
{\mu\over2}+q(0,j)={N-1\over4}-{j(N-j)\over2N}.
\nonumber
\\
&&
q(i,j)=q(0,\vert i-j\vert)=q(0,-\vert i-j\vert).
\nonumber
\end{eqnarray}
We remark that the constants $C_{\La}$ in this section are
different from those in Theorem \ref{Thomae} by $\pm1$ times
because of the reordering of the difference products.
The properties of the constants remain same.

\subsection{$N=2$}
\par
We consider the hyperelliptic curve
$s^2=\prod_{j=1}^{2m}(z-\la_l)$.
Let $\La=(\La_0,\La_1)$ with
\begin{eqnarray}
&&
\La_0=\{i_1<\cdots<i_m\},
\quad
\La_1=\{j_1<\cdots<j_m\}.
\nonumber
\end{eqnarray}
We have
\begin{eqnarray}
&&
q(0,0)={1\over8},
\quad,
q(0,1)=-{1\over8},
\quad
\mu={1\over4}.
\nonumber
\end{eqnarray}

The Thomae formula is

\begin{eqnarray}
&&
\theta[e_\La](0)^4=C_{\La}(\det A)^2
\prod_{k<l}(\la_{i_k}-\la_{i_l})(\la_{j_k}-\la_{j_l}).
\nonumber
\end{eqnarray}
This is the original Thomae formula in which case 
$C_{\La}^2=(2\pi)^{-4(m-1)}$.

\subsection{$N=3$}
Let $\La=(\La_0,\La_1,\La_2)$.
We have
\begin{eqnarray}
&&
q(0,0)={2\over9},
\quad
q(0,1)=q(0,2)=-{1\over9},
\quad
\mu={5\over9}.
\nonumber
\end{eqnarray}

Then

\begin{eqnarray}
&&
\theta[e_\La](0)^6=C_{\La}(\det A)^3
\big((\La_0\La_0)(\La_1\La_1)(\La_2\La_2)\big)^3
(\La_0\La_1)(\La_1\La_2)(\La_0\La_2).
\nonumber
\end{eqnarray}

Here if

\begin{eqnarray}
&&
\La_i=\{i_1<\cdots<i_m\},
\quad
\La_j=\{j_1<\cdots<j_m\},
\nonumber
\end{eqnarray}

then

\begin{eqnarray}
&&
(\La_i\La_i)=\prod_{k<l}(\la_{i_k}-\la_{i_l}),
\quad
(\La_i\La_j)=\prod_{k,l=1}^m(\la_{i_k}-\la_{j_l}).
\nonumber
\end{eqnarray}
Our result shows that $C_{\La}^6$ does not depend on $e_\La$.

\subsection{$N=4$}
Let $\La=(\La_0,\cdots,\La_3)$.
We have
\begin{eqnarray}
&&
q(0,0)={5\over16},
\quad
q(0,1)=q(0,3)=-{1\over16},
\quad
q(0,2)=-{3\over16},
\quad
\mu={7\over8}.
\nonumber
\end{eqnarray}

Thomae formula is
\begin{eqnarray}
\theta[e_\La](0)^8
&=&
C_{\La}(\det A)^4
\big((\La_0\La_0)(\La_1\La_1)(\La_2\La_2)(\La_3\La_3)\big)^6
\nonumber
\\
&&
\big((\La_0\La_1)(\La_1\La_2)(\La_2\La_3)(\La_0\La_3)\big)^4
\big((\La_0\La_2)(\La_1\La_3)\big)^2.
\nonumber
\end{eqnarray}
The constant $C_{\La}^{8}$ does not depend on $e_\La$.

\subsection{$N=5$}
Let $\La=(\La_0,\cdots,\La_4)$.
We have
\begin{eqnarray}
&&
q(0,0)={2\over5},
\quad
q(0,1)=q(0,4)=0,
\quad
q(0,2)=q(0,3)=-{1\over5},
\quad
\mu={6\over5}.
\nonumber
\end{eqnarray}

Thomae formula is
\begin{eqnarray}
\theta[e_\La](0)^{10}
&=&
C_{\La}(\det A)^5
\big((\La_0\La_0)(\La_1\La_1)(\La_2\La_2)(\La_3\La_3)(\La_4\La_4)\big)^{10}
\nonumber
\\
&&
\big((\La_0\La_1)(\La_1\La_2)(\La_2\La_3)(\La_3\La_4)(\La_0\La_4)\big)^6
\nonumber
\\
&&
\big((\La_0\La_2)(\La_1\La_3)(\La_2\La_4)(\La_0\La_3)(\La_1\La_4)\big)^4.
\nonumber
\end{eqnarray}
The constant $C_{\La}^{10}$ does not depend on $e_\La$.
\vskip1cm

\section{Concluding Remarks}
\par
In this paper we have given a rigorous proof of the
generalized Thomae formula for ${\bf Z}_N$ curves
which was previously discovered 
by Bershadsky and Radul \cite{BR1,BR2} in the
study of a conformal field theory.
Here let us make a comment on the related subjects.

There are several papers( \cite{G,Far} and references therein)
studying the generalization
of $\la$ function of the elliptic curves to ${\bf Z}_N$
curves by studying the cross ratios of four points
on a Riemann surface. In those approaches the only ratios
of theta constants appear and Thomae type formula is not used.
However in the Smirnov's theta formula for the solutions of
$sl_2$ Knizhnik-Zamolodchikov equation on level $0$,
Thomae formula is needed.

Our strategy to prove the generalized Thomae formula here,
which is similar to that of \cite{BR1,BR2}, is the
comparison of algebraic and analytic expressions of several
quantities.
In the hyperelliptic case this can be considered as a part of
the more general comparison of algebraic and analytic construction of
Jacobian varieties due to Mumford \cite{M}.
It will be interesting to study the integrable system
associated with ${\bf Z}_N$ curves and to study the
generalization of Thomae type formula for spectral curves.

In fact Thomae \cite{T1}, Fuchs \cite{F1} derived differential
equations satisfied by theta constants with respect to
branch points in a more general
setting and they could integrate them completely
only in the case of hyperelliptic curves.
The Thomae formula for ${\bf Z}_N$ curves provides a new example
which is integrable.

The ${\bf Z}_N$ curves and the $1/2N$ periods in the generalized
Thomae formula are related with the Lie algebra $sl_N$ and
the weight zero subspace of the tensor products of the vector
representation. Hence it is natural to expect that the Thomae
type formula has a good description in terms of Lie algebras and their
representations.

For the evaluation of the constant $C_{\La}$ we need to know
the explicit description of canonical cycles of ${\bf Z}_N$ curve.
So far we could describe a canonical basis only in the
case of $N=3$ (except $N=2$).

\vskip4mm

\noindent
{\Large\hskip4mm Acknowledgement }
\vskip4mm
\noindent
We would like to thank Koji Cho, Norio Iwase, Yasuhiko Yamada
for the illuminating discussions.

\end{document}